\documentclass[12pt,twoside,a4paper]{article}
\title{On the volume of a line bundle}
\author{S\'ebastien Boucksom}

\usepackage{latexsym}
\begin{document}
\maketitle

\newcommand{\er}{\mathbf{R}}
\newcommand{\ku}{\mathbf{Q}}
\newcommand{\ze}{\mathbf{Z}}
\newcommand{\co}{\mathbf{C}}
\newcommand{\oh}{\mathcal{O}}
\newcommand{\ih}{\mathcal{I}}
\newcommand{\ep}{\varepsilon}
\newcommand{\de}{\delta}
\newcommand{\ti}{\widetilde}
\newcommand{\ddbar}{dd^c}
\newtheorem{theo}{Theorem}[section]
\newtheorem{cor}[theo]{Corollary}
\newtheorem{defi}[theo]{Definition}
\newtheorem{prop}[theo]{Proposition}
\newtheorem{lem}[theo]{Lemma}

$Author's$ $address$: Institut Fourier, 100 rue des Maths, BP74, 38402 Saint-Martin d'H\`eres Cedex, France.\\
$e-mail$: sbouckso@ujf-grenoble.fr\\
$Abstract$: Using a result of Fujita on approximate Zariski decompositions and the singular version of Demailly's holomorphic Morse inequalities as obtained by Bonavero, we express the volume of a line bundle in terms of the absolutely continuous parts of all the positive curvature currents on it, with a way to pick an element among them which is most homogeneous with respect to the volume. This enables us to introduce the volume of any pseudoeffective class on a compact K\"ahler manifold, and Fujita's theorem is then extended to this context.  
\section{Introduction}

For a holomorphic line bundle $L$ on a compact complex manifold $X$ of dimension
$n$, one defines the $volume$ of $L$ as
$$v(L):=\limsup_{k\to+\infty}\frac{n!}{k^n}h^0(X,kL).$$
Then $L$ is big, i.e.~it has maximal Kodaira-Iitaka dimension
$\kappa(X,L)=n$, exactly when $v(L)>0$. If $L$ is an ample line bundle, the
combination of Serre's 
vanishing theorem and the asymptotic Riemann-Roch formula shows that 
$v(L)=L^n$, where $L^n$ is the $n$-fold intersection number 
$\int_Xc_1(L)^n$. When $L$ is merely a nef (numerically effective) line 
bundle and $X$ is K\"ahler, one can show using Demailly's Morse
inequalities that $h^q(X,kL)=o(k^n)$ for $q>0$, so that the
Riemann-Roch formula again yields $v(L)=L^n$ in that case. 

\noindent Recall now that a line bundle $L$ is said to have an algebraic Zariski
decomposition if there exists two $\ku$-line bundles $N$ and $E$, which are nef and
effective respectively, such that

\noindent $(i)$ $L=N+E$ as $\ku$-line bundles.

\noindent $(ii)$ $H^0(X,kL)=H^0(X,kN)$ for every positive integer $k$ 
clearing up the denominator of $N$.

\noindent
Then one of course has $v(L)=v(N)=N^n$, so that the knowledge of the nef 
part $N$ enables one to compute the volume of $L$. On a surface, such a 
decomposition exists, but this is no longer true in general, even if 
one allows $X$ to be blown up; however the following result, due to T.Fujita ([Fu94], 
cf.~also [DEL00]), shows that in some sense it is still true asymptotically.

\begin{theo}[Approximate Zariski decomposition] Let $L$ be a big line bundle 
on the projective manifold $X$. Then for every $\ep$ there exists a finite 
sequence of blowing-ups with smooth centers $\mu:\ti X\to X$ and two $\ku$-line
bundles $A$ and 
$E$ on $\ti X$ which are ample and effective respectively (the data depends 
on $\ep$) such that:

(i) $L=A+E$ as $\ku$-line bundles.

(ii) $v(A)\leq v(L)\leq v(A)+\ep$

\end{theo}
One would like to say that in general an algebraic Zariski decomposition of 
$L$ could be obtained by letting $\ep$ tend to zero in the above theorem, 
which cannot be done in an algebraic context. Our first endeavour is 
to translate the theorem of Fujita into the language of singular Hermitian 
metrics on $L$ (here the main tool is the Calabi-Yau theorem), so as
to get in the limit a ``most homogeneous'' metric on $L$ with respect to the 
volume. More precisely, we prove the following:

\begin{theo} Let $L$ be a pseudoeffective line bundle on a compact K\"ahler manifold $X$. Then:

\noindent $(i)$ $v(L)=\sup_T\int_XT_{ac}^n$ for $T$ ranging among the closed 
positive $(1,1)$-currents in the cohomology class $c_1(L)$.

\noindent $(ii)$ Given a K\"ahler form $\omega$ on $X$, normalized so that 
$\int_X\omega^n=1$, there exists a closed positive $(1,1)$-current $T\in 
c_1(L)$ such that
$T_{ac}^n(x)=v(L)\omega^n(x)$ for almost every $x\in X$.
\end{theo}
Recall that a line bundle $L$ is said to be pseudoeffective if it can be 
equipped with a singular Hermitian metric $h$ with positive curvature 
current $\Theta_h(L)$, or equivalently if the (de Rham) cohomology class 
$c_1(L)$ contains a closed positive $(1,1)$-current $T$. $T_{ac}$ is then  
the absolutely continuous part of $T$ in its Lebesgue decomposition 
$T=T_{ac}+T_{sg}$. $T_{ac}$ is considered as a $(1,1)$-form with $L^1_{loc}$ 
coefficients, so that the product $T_{ac}^n$ is meant pointwise.

In particular, Theorem 1.2 shows that the volume $v(L)$ of $L$ only depends on 
the first Chern class $c_1(L)$. It therefore seems natural to introduce the 
following

\begin{defi} Let $a$ be a pseudoeffective $(1,1)$-cohomology class on $X$ 
compact K\"ahler. Then one defines its volume by

$$v(a):=\sup_T\int_XT_{ac}^n$$
for $T$ ranging over the closed positive $(1,1)$-currents in $a$.
\end{defi}
We show that the volume $v(a)$ is always finite, that it is continuous on 
the whole pseudoeffective cone and that Theorem 1.2 still holds in this more 
general context. Furthermore, we give the corresponding version of Fujita's 
theorem:

\begin{theo} Let $X$ be a compact K\"ahler manifold, and let $a$ be a 
pseudoeffective $(1,1)$-cohomology class on $X$ with $v(a)>0$. Then for 
every $\ep>0$ there exists a finite sequence of blowing-ups with smooth centers $\mu:\ti X\to 
X$, a K\"ahler class $\omega$ and an effective real divisor $D$ on $\ti X$ 
such that

\noindent $(i)$ $\mu^{\star}a=\omega+D$ as cohomology classes.

\noindent $(ii)$ $v(\omega)=\omega^n\leq v(a)\leq v(\omega)+\ep$
\end{theo}
In particular, a pseudoeffective class with positive volume is $big$, in the 
sense that it contains a K\"ahler current (cf.~definitions below).\\

\section{Preliminary tools.}

\subsection{Currents with analytic singularities}

In this section, $X$ denotes a complex $n$-dimensional manifold.
\noindent
Recall first that

\noindent $(i)$ A closed real $(1,1)$-current $T$ on $X$ is said to be 
almost positive if some continuous real $(1,1)$-form $\gamma$ can be found
such that $T\geq\gamma$, and that a function $\varphi$ in $L^1_{loc}(X)$ is 
almost plurisubharmonic (psh for short) if its Hessian $dd^c\varphi$ is an almost positive 
current (here $dd^c=\frac{i}{\pi}\partial\overline{\partial}$). This latter 
property is equivalent to the fact that $\varphi$ can locally be written as 
a sum of a plurisubharmonic function and a smooth one.

\noindent $(ii)$ A closed $(1,1)$-current $T$ is called a K\"ahler current 
if one has $T\geq\omega$ for some Hermitian form $\omega$ on $X$. For 
example, a line bundle $L$ on $X$ complex compact is big if and only if it 
admits a curvature current which is K\"ahler. That the condition is 
sufficient is a consequence of a result of Bonavero [Bon93] which will be 
stated below (cf.~also [JS93], or [Dem96] for the easier case where $X$ is
K\"ahler). In the other direction, choose $k$ such that $H^0(X,kL)$ 
generates $1$-jets generically, and fix some smooth metric $h_{\infty}$ on 
$L$. Then the metric $h:=h_{\infty}e^{-2\varphi}$ has K\"ahler curvature 
current if one takes\\ 
$\varphi:=\frac{1}{2k}\log(\sum|\sigma_j|_{h_{\infty}}^2)$ for some basis 
$(\sigma_j)$ of $H^0(X,kL)$. But the fact that this metric is obtained 
algebraically has a nice consequence on the singularities of $h$. In fact, 
it has algebraic singularities in the following sense:

\begin{defi}
$(i)$ Given a coherent ideal sheaf $\ih$ and $c>0$, we say that a function 
$\varphi$ has singularities of type $(\ih,c)$ if it can locally be written 
as $\varphi=\frac{c}{2}\log(\sum|f_j|^2)+\theta$ for some local 
generators $(f_j)$ of $\ih$ and some smooth function $\theta$.

\noindent $(ii)$ We say that $\varphi$ has analytic (resp.~algebraic) 
singularities if it has singularities of type $(\ih,c)$ for some ideal sheaf 
$\ih$ and some real (resp.~rational) $c>0$.
\end{defi}
A function with analytic singularities is automatically almost psh, and we 
can use the same terminology about the singularities of a closed almost 
positive current by looking at its local potentials. One can translate 
algebraic data into potential theoretic ones by means of the

\begin{prop} Given $c>0$ and a coherent ideal sheaf $\ih$, there exists an 
almost psh function $\varphi$ on $X$ with singularities of type $(\ih,c)$.
\end{prop}
The proof consists merely in gluing together by means of a partition of 
unity functions of the type $\frac{c}{2}\log(\sum_j|g_j|^2)$, where the $g_j$ are local generators of $\ih$. To do this, cover $X$ by a finite number of open 
sets $U_k$ such that $\ih$ is generated by $f_{k,j}$, $j=1,...,N_k$ on 
$U_k$. Then take $\theta_k$ smooth, positive, with support in $U_k$, and 
such that $\sum_k\theta_k^2=1$. Then one sets 
$\varphi=\frac{c}{2}\log(\sum_k\theta_k^2 e^{\varphi_k})$ with 
$\varphi_k=\log(\sum_j|f_{k,j}|^2)$.

\noindent
Let us stress that the data $(\ih,c)$ is not uniquely determined for a function
$\varphi$ with analytic singularities (for instance $\ih$ can always
be replaced by its integral closure). We have the following obvious functorial
property of functions with analytic 
singularities:

\begin{prop} Let $f:X\to Y$ be a holomorphic map between complex \\manifolds, 
and let $\varphi$ have singularities of type $(\ih,c)$. Then 
$f^{\star}\varphi$ has singularities of type $(f^{-1}\ih,c)$.
\end{prop}

\begin{cor}[Resolution of singularities] Given $\varphi$ with singularities 
described by $(\ih,c)$ on $X$ compact complex, there exists a finite 
sequence of blowing-ups $\mu:\ti X\to X$ such that $\mu^{\star}\varphi$ has 
singularities along a divisor of $\ti X$ with normal crossings.
\end{cor}
This follows by first blowing up $X$ along $\ih$ and then applying 
Hironaka's resolution of singularities.

\noindent Now recall Siu's decomposition formula: given a closed almost 
positive current $T$ of bidimension $(p,p)$, one can write 
$$T=\sum_j\nu_j[A_j]+R$$where the $A_j$ are irreducible $p$-dimensional 
analytic subsets, 

\noindent $\nu_j=\min_{x\in A_j}\nu(T,x)$ is the generic Lelong 
number of $T$ along $A_j$, and $R$ is a closed almost positive 
$(1,1)$-current such that $\dim E_c(R)<p$ for every $c>0$ 

\noindent ($E_c(R)=\{x\in 
X,\nu(R,x)\geq c\}$ is the Lelong sublevel set, which is analytic by a
well known result of Y.T.Siu.)
Then we claim the following

\begin{prop} For a closed $(1,1)$-current with singularities of type 
$(\ih,c)$, the Siu decomposition coincides with the Lebesgue decomposition, 
i.e.~$T_{ac}=R$ and $T_{sg}=\sum_j\nu_j[A_j]$. Furthermore, $T_{sg}$ is a 
finite sum which corresponds to the divisorial part of $V(\ih)$.
\end{prop}
Indeed, working locally, suppose that 
$\varphi=\frac{1}{2}\log(|f_1|^2+...+|f_N|^2)$, and let $g$ be the g.c.d. of 
the $f_j$'s. Then $dd^c\varphi=D+R$ where $D=dd^c\log|g|$ is the integration 
current along the divisor $\{g=0\}$ (by the Poincar\'e-Lelong formula) and 
$R=\frac{1}{2}dd^c\log(|g_1|^2+...+|g_N|^2)$, thus it suffices to show that 
$R$ is absolutely continuous. But it is even smooth outside the set $A$ of 
common zeroes of the $g_j$'s, and carries no mass on $A$ since codim$A\geq 
2$ by construction, thus it has $L^1_{loc}$ coefficients.

\noindent  Note that as a consequence $T_{ac}$ and $T_{sg}$ are closed for 
$T$ with analytic singularities, which is not true in general.

\noindent  Now in general one has the following result due to Demailly 
[Dem92]:

\begin{theo}[Approximation with algebraic singularities]
Let $T$ be a closed almost positive $(1,1)$-current on $X$ compact complex equipped with some Hermitian form $\omega$. 
Write $T=\alpha+dd^c\varphi$ for $\alpha$ smooth and $\varphi$ almost psh, 
and suppose also given some continuous real form $\gamma$ such that 
$T\geq\gamma$. Then there exists a sequence $\varphi_k$ of almost psh 
functions such that:

\noindent $(i)$ $\varphi_k$ has algebraic singularities.

\noindent $(ii)$ $\varphi_k$ decreases to $\varphi$ everywhere.

\noindent $(iii)$ The Lelong numbers $\nu(\varphi_k,x)$ converge to 
$\nu(\varphi,x)$ uniformly with respect to $x$.

\noindent $(iv)$ Setting $T_k=\alpha+dd^c\varphi_k$, one has 
$T_k\geq\gamma-\ep_k\omega$ for some
sequence $\ep_k$ decreasing to $0$.
\end{theo}
For our concerns, we need a slight improvement of this result, in order to 
take care of the absolutely continuous parts in the convergence.

\begin{cor} In the above proposition, we can also impose

\noindent $(v)$ $T_{k,ac}\to T_{ac}$ almost everywhere.
\end{cor}
Proof: We appeal to a different approximation result [Dem82] taking care of 
the absolutely continuous parts (it is essentially a convolution): there 
exists a sequence $\psi_j$ of smooth functions on $X$ such that, setting 
$S_j=\alpha+dd^c\psi_j$, we have:

\noindent $(a)$ $\psi_j$ decreases to $\varphi$ everywhere.

\noindent $(b)$ $S_j\geq\gamma-C\lambda_j\omega$ for some $C>0$ and 
continuous functions $\lambda_j$ decreasing pointwise to the Lelong number,

\noindent $(c)$ $S_j(x)\to T_{ac}(x)$ almost everywhere.

What we will do is glue the two sequences $\varphi_k$ and $\psi_j$ in 
the following fashion: denote by $A_k$ the polar set of $\varphi_k$ and choose an 
arbitrary sequence $C_k$ of positive reals increasing to $+\infty$, and a 
sequence of rationals $\de_k$ decreasing to 0.
Now observe that $U_k:=\{\varphi_k<-(C_k+1)/\de_k\}$ is an open 
neighbourhood of $A_k$ such that we have
$\varphi_k<(1-\delta_k)\varphi_k-C_k-1/2$ on $\overline{U}_k$, so that $$\varphi\leq\varphi_k<(1-\de_k)\varphi_k-C_k-1/2 \mbox{ on
  the compact }\overline{U}_k.$$
Thus $(a)$ and the continuity of $\psi_j$ show that for $j_k$ big enough, 
we have \\$\psi_{j_k}<(1-\de_k)\varphi_k-C_k-1/2$ on
$\overline{U}_k$. Now we select a smaller open neighbourhood
$W_k\subset\subset U_k$ of $A_k$, and we set:
$$\theta_k:=\cases{(1-\de_k)\varphi_k-C_k& on $W_k$,\cr
\max_{\eta}((1-\de_k)\varphi_k-C_k,\psi_{j_k})& on $X-W_k$\cr}$$ 
with $\max_{\eta}(x,y)$ denoting a regularized maximum function.

\noindent
If $\eta$ is chosen so small that $\max_{\eta}(x,y)=x$ when $y<x-1/2$ then 
the two parts to be glued coincide on some neighbourhood of $\partial W_k$, 
so that $\theta_k$ is almost plurisubharmonic with algebraic singularities 
on $X$.

\noindent The gluing property of plurisubharmonic functions shows that 
$\alpha+dd^c\theta_k$ will have a lower bound going to 0 for $k\to+\infty$ 
if we can show that this is the case for $\alpha+dd^c\psi_{j_k}$ on $X-W_k$. 
But we have $\nu(\varphi_k,x)=0$ for $x$ in this set, thus $(iii)$ and $(b)$ 
together give the result.

Now we have to see that $(\alpha+dd^c\theta_k)_{ac}(x)\to T_{ac}(x)$ almost 
everywhere. To do this, notice that if $\varphi(x)>-\infty$, then $x$ cannot 
be in $U_k$ for $k$ big enough, since otherwise $\varphi_k(x)\leq-(C_k+1)$ 
for $k$ big enough, which would yield $\varphi(x)= -\infty$ by $(ii)$. 
Furthermore we have $$(1-\de_k)\varphi_k(x)-C_k<\varphi(x)-1/2\leq\psi_j-1/2$$
for $k$ big enough (since $C_k\to+\infty$) and $j$ arbitrary (by $(a)$), thus in
particular for 
$j=j_k$. Continuity then gives $(1-\de_k)\varphi_k-C_k<\psi_{j_k}-1/2$ on some neighbourhood of $x$ (depending on $k$) contained in $X-W_k$,
so that $\theta_k=\psi_{j_k}$ on 
this neighbourhood. Thus for every $x$ outside the polar set of $\varphi$ 
(which has measure 0) we have $(\alpha+dd^c\theta_k)(x)=S_{j_k}(x)$ for $k$ big
enough, and the result follows by $(c)$.

To conclude this section, let us illustrate the above ideas by reproducing 
the following result from [DP01]:

\begin{prop} A compact complex manifold $X$ is in the Fujiki class
  $\cal C$ if, 
and only if, it carries a K\"ahler current.
\end{prop}
Recall that $X$ is said to be in the Fujiki class $\cal C$ if it can be modified
into 
a K\"ahler manifold. Since the push-forward of a K\"ahler form by a 
modification is a K\"ahler current, one direction is clear. Assume now that 
the compact complex manifold $X$ carries a K\"ahler current $T\geq\omega$ 
for some Hermitian form $\omega$. Using the above approximation result, one 
can assume that $T$ has analytic singularities. By corollary 2.4, we get a 
finite sequence of blowing-ups $\mu:\ti X\to X$ such that $\mu^{\star}T$ 
decomposes into a sum $\alpha+D$, with $\alpha$ a smooth form and $D$ a real 
effective divisor. Then $\alpha\geq\mu^{\star}\omega$ is semi-positive, but 
not definite positive everywhere $a$ $priori$. Here one uses the 
following

\begin{lem} Let $f:X\to Y$ be the blow-up of $Y$ along a submanifold $V$, let 
$\omega$ be a Hermitian form on $X$ and let $E$ denote the exceptional divisor 
of $f$. Then there exists a smooth closed form $u$ cohomologous to $E$ such 
that $\mu^{\star}\omega-\ep u$ is definite positive everywhere for $\ep>0$ 
small enough.
\end{lem}
To see this, recall that $\oh_E(-E)$ is isomorphic to $\oh_{P(N_{V/Y})}(1)$ 
on $P(N_{V/Y})=E$, hence is relatively ample over $V$, so that if we choose 
some smooth form $\beta$ on $X$ cohomologous to $[E]$, there exists a smooth 
function $\varphi$ on $E$ such that $-(\beta_{|E}+dd^c\varphi)$ is definite 
positive along the fibers of $E\to V$. Now we extend $\varphi$ smoothly to 
the whole of $X$, and set $u:=\beta+dd^c\varphi$. If we replace $\varphi$ 
by $\varphi-Cd(x,E)$ for $C\gg0$, we can even arrange so that $-u$ be 
positive definite on $N_{E/X}$.
Now we note that $\mu^{\star}\omega$ is semi-positive everywhere and 
strictly positive on the ``horizontal'' directions of $E\to V$, so that 
$\mu^{\star}\omega-\ep u$ will be positive definite on a neighbourhood $E$ for 
$\ep>0$ small enough. There remains to choose $\ep>0$ even smaller to ensure 
that $\mu^{\star}\omega-\ep u$ be positive outside this neighbourhood, which 
can be done because $\mu^{\star}\omega$ has a uniform strictly positive 
lower bound there.

A repeated application of this lemma will therefore yield a 
smooth closed form $v$ cohomologous to an effective linear combination with small coefficients of the exceptional divisors of $\mu$ such 
that $\alpha-v$ is definite positive, hence a K\"ahler form, for
$\ep>0$ small enough.

\subsection{Boundedness of the mass.}

Here we are interested in the control of $\int_XT_{ac}^n$ for a closed 
positive $(1,1)$-current. It is by no means clear that such an integral is 
convergent in general, and we show that things go well in the K\"ahler case.

\begin{lem} Let $T$ be any closed positive $(1,1)$-current on $X$ compact 
complex. Then the Lelong numbers $\nu(T,x)$ of $T$ can be bounded by a 
constant depending only on the $dd^c$-cohomology class of $T$.
\end{lem}
Indeed for $\omega$ a Gauduchon form on $X$ (i.e.~a Hermitian form with 
$dd^c(\omega^{n-1})=0$), one has by definition that $\nu(T,x)$ is (up to a 
constant depending on $\omega$ near $x$) the limit for $r\to 0_+$ of
$$\nu(T,x,r):=\frac{(n-1)!}{(\pi
r^2)^{n-1}}\int_{B(x,r)}T\wedge\omega^{n-1},$$known to be an increasing 
function of $r$. Thus if we choose $r_0$ small enough to ensure that each 
ball $B(x,r_0)$ is contained in a coordinate chart, we get 
$\nu(T,x)\leq\nu(T,x,r_0)\leq C\int_XT\wedge\omega^{n-1}$, a quantity 
depending only on the $dd^c$ cohomology class $\{T\}$.

\begin{lem}[Uniform boundedness of the mass)]Let $T$ be a closed positive 
$(1,1)$-current on $(X,\omega)$ compact K\"ahler. Then the integrals 
$\int_XT_{ac}^k\wedge\omega^{n-k}$ are finite for each $k=0,...,n$ and can 
be bounded in terms of $\omega$ and the cohomology class of $T$ only.
\end{lem}
Proof: By [Dem82], there exists a sequence $T_j$ of smooth $(1,1)$-forms in
the same class as $T$ and a constant $C>0$ depending only on $(X,\omega)$ 
such that:

\noindent $(i)$ $T_j\to T$ weakly,

\noindent $(ii)$ $T_j(x)\to T_{ac}(x)$ almost everywhere, and

\noindent $(iii)$ $T_j\geq -C\lambda_j\omega$ for some continuous functions 
$\lambda_j$
such that $\lambda_j(x)$ decreases to $\nu(T,x)$ when $j\to+\infty$, for 
every $x$.

By using the preceding lemma, one can thus find a
constant also denoted by $C$ and depending on $\omega$ and the cohomology class 
$\{T\}$ only such that $T_j+C\omega\geq 0$. But now
$$\int_X(T_j+C\omega)^k\wedge\omega^{n-k}=\{T+C\omega\}^k\{\omega\}^{n-k}$$
does not depend on $j$, so the result follows by Fatou's lemma.

Here, even though the approximation result of Demailly does not need
the K\"ahler assumption, the argument definitely fails in its last part 
without it. However it still holds true if $\dim X\leq 2$, as one sees
in the above proof.

\section{Volume of a line bundle}

\subsection{A Morse-type inequality}
In this part, we will prove the following:

\begin{prop} For $L$ pseudoeffective on $X$ compact K\"ahler and $T\in 
c_1(L)$ positive one has:
$$v(L)\geq\int_XT_{ac}^n$$
\end{prop}
In order to get this, we will appeal to the singular holomorphic Morse 
inequalities, but to state this result we first need some terminology.

\noindent $(i)$ The Nadel multiplier ideal sheaf $\ih(T)$ of an almost 
positive closed $(1,1)$-current $T$ is defined as the sheaf of germs of 
holomorphic functions $f$ such that $|f|^2e^{-2\varphi}$ is locally 
integrable, for some (hence any) $\varphi\in L^1_{loc}$ such $T=dd^c\varphi$ 
locally. This sheaf is coherent, as is well known (cf.~e.g.~[Dem96]).

\noindent $(ii)$ The $q$-index set of an almost positive closed 
$(1,1)$-current $T$ is the set of $x\in X$ such that the absolutely 
continuous part $T_{ac}$ of $T$ has exactly $q$ negative eigenvalues at $x$. 
This set is denoted by $X(T,q)$, and we also write $X(T,\leq q)$ for the 
union of the $X(T,j)$'s for $j=0,...q$. These sets are only defined up to a 
null measure set, but since we shall only integrate absolutely 
continuous forms on them, it is not really annoying here.

Now we can state the following result, due to Bonavero [Bon93]:
\begin{theo}[Singular Morse inequalities] Let $L$ be any holomorphic line 
bundle on the compact complex $n$-dimensional manifold $X$, and $T\in 
c_1(L)$ be some closed $(1,1)$-current with algebraic singularities. Then 
the following holds:
$$h^0(X,\oh(kL)\otimes\ih(kT))-h^1(X,\oh(kL)\otimes\ih(kT))\geq\frac{k^n}{n!}\int_{X(T,\leq
1)}T_{ac}^n-o(k^n)$$
for $k\to+\infty$.
\end{theo}
Since one has
$$h^0(X,\oh(kL))\geq h^0(X,\oh(kL)\otimes\ih(kT))$$ $$\geq
h^0(X,\oh(kL)\otimes\ih(kT))-h^1(X,\oh(kL)\otimes\ih(kT))$$Theorem 3.2 implies:
\begin{cor} For any line bundle $L$ on the compact complex $X$, one has
$$v(L)\geq\int_{X(T,\leq 1)}T_{ac}^n$$for every $T\in c_1(L)$ with algebraic 
singularities.
\end{cor}
We can now conclude the proof of Proposition 3.1. In fact, we choose a 
sequence $T_k$ of currents with algebraic singularities as in Corollary 2.7 
(here $\gamma=0$), and we denote by $\lambda_1\leq ...\leq\lambda_n$ (resp.~
$\lambda^{(k)}_1\leq ...\leq\lambda_n^{(k)}$) the eigenvalues of $T_{ac}$ 
(resp.~$T_{k,ac}$) with respect to $\omega$. We have by assumption that 
$\lambda_1\geq 0$, $\lambda^{(k)}_1\geq-\ep_k$ and 
$\lambda^{(k)}_j(x)\to\lambda_j(x)$ almost everywhere. We can certainly 
assume that $\int_XT_{ac}^n>0$, which means that the set 
$A:=\{\lambda_1>0\}$ has positive measure. For each small $\delta>0$, 
Egoroff's lemma gives us some $B_{\delta}\subset A$ such that 
$\lambda^{(k)}_1\to\lambda_1$ uniformly on $B_{\delta}$ and also 
$A-B_{\delta}$ has measure less than $\delta$. Thus we see that 
$B_{\delta}\subset X(T_k,0)$ for $k$ big enough, and consequently

\noindent
$\limsup\int_{X(T_k,0)}T_{k,ac}^n\geq\int_{B_{\delta}}\liminf 
T_{k,ac}^n=\int_{B_{\delta}}T_{ac}^n$, using Fatou's lemma.
Letting now $\delta$ tend to 0, we get
$$\limsup\int_{X(T_k,0)}T_{k,ac}^n\geq\int_AT_{ac}^n=\int_XT_{ac}^n.$$ Since 
by Corollary 3.3 above we have $\int_{X(T_k,0)}T_{k,ac}^n\leq 
v(L)-\int_{X(T_k,1)}T_{k,ac}^n$ for every $k$, the proof of the inequality 
will be over if we can show that\\$-\int_{X(T_k,1)}T_{k,ac}^n\to 0$. But we
observe the following inequalities on $X(T_k,1)$:

$$0\leq -T_{k,ac}\leq n\ep_k\omega\wedge(T_{k,ac}+\ep_k\omega)^{n-1},$$
from which we get

$$0\leq-\int_{X(T_k,1)}T_{k,ac}^n\leq 
n\ep_k\int_X\omega\wedge(T_{k,ac}+\ep_k\omega)^{n-1}.$$ Now the last 
integral is bounded uniformly in terms of $\{T\}$ and $\omega$ only by the 
``boundedness of the mass'' (Lemma 2.11), which ends the proof.

It is worth noting that the K\"ahler assumption is needed precisely for this
last lemma. Consequently, Proposition 3.1 is also true on any surface.

\subsection{The theorem of Calabi-Yau.}
We fix a K\"ahler form $\omega$ with $\int_X\omega^n=1$. We seek a 
positive $T$ in $c_1(L)$ such that $T_{ac}^n\geq v(L)\omega^n$ $(\star)$. If 
this is done, the proof of Theorem 1.2 is over, since the fact that 
$\int_XT_{ac}^n\leq\int_Xv(L)\omega^n$ (by Proposition 3.1) implies that
$(\star)$ is an equality a.e. First let us recall the fundamental
result proved in [Yau78]:

\begin{theo} [Aubin-Calabi-Yau] Let $(X,\omega)$ be a compact K\"ahler 
manifold, and assume that $\int_X\omega^n=1$. Then given a K\"ahler 
cohomology class $a$, there exists a K\"ahler form $\alpha\in a$ such that
$$\alpha^n=(\int_Xa^n)\omega^n$$
pointwise.
\end{theo}
We now appeal to Fujita's Theorem 1.1, and use the notations in it. Notice 
first that if our $L$ is not big, we will have $0=v(L)\geq\int_XT_{ac}^n\geq 
0$ for every positive $T\in c_1(L)$ because of what we have already seen, so that any such 
$T$ would do. We thus assume that $L$ is big, which implies automatically 
the projectivity of $X$ since the latter will be K\"ahler and Moishezon at 
the same time.

\noindent Now we choose any K\"ahler form $\omega_{\ep}$ on $X_{\ep}$, and 
apply the Calabi-Yau theorem to $A_{\ep}$ with respect to the K\"ahler form 
$\mu_{\ep}^{\star}\omega+\delta\omega_{\ep}$, normalized adequately, and we 
get a K\"ahler form $\alpha_{\ep,\delta}\in c_1(A_{\ep})$ such that 
$$\alpha_{\ep,\delta}^n=\frac{v(A_{\ep})}{\int_{X_{\ep}}(\mu_{\ep}^{\star}\omega
+\delta\omega_{\ep})^n}(\mu_{\ep}^{\star}\omega+\delta\omega_{\ep})^n.$$Now 
consider 
$$T_{\ep,\delta}:=(\mu_{\ep})_{\star}(\alpha_{\ep,\delta}+[E_{\ep}]).$$This 
is a positive current lying in $c_1(L)$, and so standard compactness 
properties of currents show that for each $\ep$ there exists some sequence 
$\delta_k(\ep)$ decreasing to 0 such that $T_{\ep,\delta_k}$ converges to some 
positive $T_{\ep}\in c_1(L)$. For the same reason there exists a sequence 
$\ep_k$ decreasing to 0 such that $T_{\ep_k}$ converges to some positive 
$T\in c_1(L)$. We will show that this $T$ satisfies $T_{ac}^n\geq 
v(L)\omega^n$ almost everywhere, which will conclude the proof. This will be 
handled by the following semicontinuity property:

\begin{lem} Let $T_k$ be a sequence of positive $(1,1)$-currents converging 
weakly to $T$. Then one has $T_{ac}^n\geq\limsup T_{k,ac}^n$ a.e.
\end{lem}
Proof: We fix a $\omega$ a Hermitian form, and for $\alpha$ a positive 
$(1,1)$-form we denote by $\det(\alpha)$ the determinant of $\alpha$ with 
respect to $\omega$, that is $\det(\alpha)=\alpha^n/\omega^n$. Since the 
result is local, we may consider a regularizing sequence $(\rho_j)$. The 
concavity of the function $A\to\det(A)^{1/n}$ on the convex cone of Hermitian 
semi-positive matrices of size $n$ then yields the second of the following 
inequalities
$$\det(T_k\star\rho_j)^{1/n}\geq \det(T_{k,ac}\star\rho_j)^{1/n}\geq 
\det(T_{k,ac})^{1/n}\star\rho_j.$$Since convolution transforms a weak 
convergence into a $C^{\infty}$ one, Fatou's lemma therefore implies:
$$\det(T\star\rho_j)^{1/n}\geq(\liminf 
\det(T_{k,ac})^{1/n})\star\rho_j.$$Now Lebesgue's theorem implies that 
$T\star\rho_j\to T_{ac}$ a.e., thus we get
$$\det(T_{ac})^{1/n}\geq\liminf \det(T_{k,ac})^{1/n},$$. We can eventually turn the $\liminf$ into a $\limsup$ by choosing subsequences pointwise.

\noindent $Remark$: this type of arguments (Calabi-Yau+convolution and
concavity) 
can already be found in [Dem93].

\section{Volume of a pseudoeffective class.}

\subsection{General properties.}

In this section, $X$ is a compact K\"ahler manifold unless otherwise 
specified.We propose to extend some results related to the volume of a line 
bundle, that is to say of an entire (or even rational) pseudoeffective 
class, to the more general case of an arbitrary pseudoeffective 
$(1,1)$-class $a$. Note that the volume $v(a)$ of $a$, as defined in the 
introduction, is always finite by Lemma 2.11 (boundedness of the mass).

\noindent First of all, let us state the following 
\begin{prop} If $f:X\to Y$ is a generically finite holomorphic map between 
compact K\"ahler manifolds and $a$ $($resp.~$b$ $)$ is a pseudoeffective 
$(1,1)$-class on $X$ $($resp.~$Y)$, then one has $v(f_{\star}a)\geq v(a)$ $($resp.
$v(f^{\star}b)=(\deg f) v(b))$.
\end{prop}
This is consequence of the following easy facts:

\begin{lem} Let $f$ be as above. Then the following properties hold:

\noindent $(i)$ If $\alpha$ is a form with $L^1_{loc}$ coefficients on $X$, 
then $f_{\star}\alpha$ (taken pointwise) is $L^1_{loc}$ on $Y$ (this only 
requires $f$ to be proper and surjective), and 
$f_{\star}(\alpha^k)=(f_{\star}\alpha)^k$ (taken pointwise) whenever 
$\alpha^k$ is also $L^1_{loc}$,

\noindent $(ii)$ If $\omega$ is an integrable top degree form on $X$ 
$($resp.~on $Y)$, then $\int_Yf_{\star}\omega=\int_X\omega$ 
$($resp.~$\int_Xf^{\star}\omega=(\deg f)\int_Y\omega)$,

\noindent $(iii)$ If $T$ is a positive $(1,1)$-current on $X$ $($resp.~on 
$Y)$, then $(f_{\star}T)_{ac}=f_{\star}(T_{ac})$ a.e. $($resp.~
$(f^{\star}T)_{ac}=f^{\star}(T_{ac})$ a.e.$)$ where $f_{\star}(T_{ac})$ 
$($resp.~$f^{\star}(T_{ac}))$ are meant pointwise.
\end{lem}

\noindent We have already seen that for a nef line bundle $L$ one has 
$v(L)=L^n$. This remains true for an arbitrary nef class:

\begin{prop} Let $a$ be a nef class on $X$. Then one has
$$v(a)=\int_Xa^n$$
\end{prop}
The proof is in two steps. First we give the following more precise nef 
version of the boundedness-of-the-mass lemma, due to Mourougane [Mou98]:

\begin{lem} Given a nef class $a$ on the K\"ahler manifold $(X,\omega)$ (or any
compact complex surface), one has for every positive $T\in a$ the following
inequalities:
$$\int_XT_{ac}^k\wedge\omega^{n-k}\leq \int_Xa^k\wedge\omega^{n-k}$$
for $k=0,...,n$.
\end{lem}
For the proof, we write as before $T=\alpha+dd^c\varphi$ with $\alpha$ 
smooth in $a$, and consider a smooth sequence $T_j=\alpha+dd^c\psi_j$ in $a$ 
such that

\noindent $(i)$ $\psi_j$ decreases to $\varphi$,

\noindent $(ii)$ $T_j\geq -C\lambda_j\omega$ for $C>0$ and $\lambda_j$ 
continuous decreasing to the Lelong number,

\noindent $(iii)$ $(T_j)(x)\to T_{ac}(x)$ almost everywhere.

\noindent Since $a$ is nef, we also have by definition a sequence of smooth
forms 
$\alpha_k=\alpha+dd^c\varphi_k$ such that $\alpha_k\geq -\ep_k\omega$  for a 
sequence $\ep_k$ decreasing to 0.
We will again glue these two sequences together to get a smooth sequence in 
$a$ with small negative part  and a good behaviour with respect to the 
absolutely continuous part, and then the same argument as for the 
boundedness of the mass will allow us to conclude.

In fact, we take a sequence of positive reals $C_k$ converging to $+\infty$, 
and we set $\theta_{k,j}:=\max_{\eta}(\varphi_k-C_k,\psi_j)$. We will show 
that for each $k$ we can choose $j_k$ big enough so that 
$\theta_k:=\theta_{k,j_k}$ has a Hessian admitting a lower bound going to 0. 
By $(ii)$, it is enough to show that given $\de>0$ and $k$, we have 
$\psi_j<\varphi_k-(C_k+1)$ near $E_{\de}(T)$ for $j$ big enough. But since 
$\varphi$ has poles on the compact set $E_{\de}(T)$, it is just a 
consequence of $(i)$ and the continuity of $\psi_j$.

As a consequence of this lemma, we get of course that $v(a)\leq\int_Xa^n$ 
for $a$ nef.
To get the converse inequality, we will prove the following result analogous 
to Theorem 1.2:

\begin{prop}[Calabi-Yau, nef case] Given a K\"ahler form $\omega$ with\\ 
$\int_X\omega^n=1$, there exists a positive $T\in a$ such that
$$T_{ac}^n=(\int_Xa^n)\omega^n$$ almost everywhere.
\end{prop}
This is of course a consequence of the Calabi-Yau theorem as follows: for 
every $\ep>0$ there exists a K\"ahler form $\alpha_{\ep}\in a+\ep\omega$ such
that 
$$\alpha_{\ep}=(\int_X(a+\ep\omega)^n)\omega^n$$pointwise.
Since the mass 
$\int_X\alpha_{\ep}\wedge\omega^{n-1}=\{a+\ep\omega\}\{\omega\}^{n-1}$ is 
bounded, there exists some weak limit $T\in a$ of $\alpha_{\ep}$. Then Lemma 3.5
yields
$T_{ac}^n\geq(\int_Xa^n)\omega^n$.
By the reverse inequality proven above, this is enough to conclude.

\subsection{A degenerate Calabi-Yau theorem.}

In this section, we prove the following more general version of Theorem 1.2:
\begin{theo}[Calabi-Yau, pseudoeffective case] If $a$ is a pseudoeffective 
class and $\omega$ a K\"ahler form with $\int_X\omega^n=1$, then there is a 
positive $T\in a$ such that
$$T_{ac}^n=v(a)\omega^n$$ almost everywhere.
\end{theo}
Proof: choose a positive $T\in a$ with $\int_XT_{ac}^n$ close to $v(a)$, and 
then a sequence $T_{\ep}\in a$ with analytic singularities such that 
$T_{\ep}\geq-\ep\omega$ and $T_{\ep,ac}\to T_{ac}$ almost everywhere. 
Let then $\mu_{\ep}:X_{\ep}\to X$ be a modification such that 
$\mu_{\ep}^{\star}T_{\ep}$ decomposes as a sum $\alpha_{\ep}+D_{\ep}$ with 
$\alpha_{\ep}$ smooth and $D_{\ep}$ a real effective divisor. Then 
$\alpha_{\ep}+\ep\mu_{\ep}^{\star}\omega$ defines a nef class on $X_{\ep}$, 
so by the nef case of Calabi-Yau, if we select an arbitrary K\"ahler form 
$\omega_{\ep}$ on $X_{\ep}$ and $\de>0$, there is a positive $T_{\ep,\de}$ 
cohomologuous to $\alpha_{\ep}+\ep\mu_{\ep}^{\star}\omega$ such that
$$T_{\ep,\de,ac}^n=\frac{v(\alpha_{\ep}+\ep\mu_{\ep}^{\star}\omega)}{v(\mu_{\ep}
^{\star}\omega+\de\omega_{\ep})}(\mu_{\ep}^{\star}\omega+\de\omega_{\ep})^n$$
almost everywhere.
Now we look at the current
$$S_{\ep,\de}:=\mu_{\ep,\star}(T_{\ep,\de}+D_{\ep})$$
on $X$. Its cohomology class is $a+\ep\omega$, so after extracting sequences 
decreasing to 0 for $\de$, $\ep$ in that order we get some current
$$S=\lim_{\ep\to 0}\lim_{\de\to 0}S_{\ep,\de}$$
in the class $a$. Since $\mu_{\ep}$ is a local isomorphism almost everywhere, we have 
$$S_{\ep,\de,ac}^n=\frac{v(\alpha_{\ep}+\ep\mu_{\ep}^{\star}\omega)}{v(\mu_{\ep}
^{\star}\omega+\de\omega_{\ep})}(\omega+\de\mu_{\ep,\star}\omega_{\ep})^n$$almost
everywhere, and Lemma 3.5 therefore implies:
$$S_{ac}^n\geq(\lim_{\ep\to 0}\int_X(T_{\ep,ac}+\ep\omega)^n)\omega^n.$$
The inner limit is greater than $\int_XT_{ac}^n$ by Fatou's lemma. Now we 
let $\int_XT_{ac}^n$ tend to $v(a)$, and choose $S_{\infty}$ some 
accumulation point of the currents $S$ obtained above, so that eventually we 
obtain $S_{\infty,ac}^n\geq v(a)\omega^n$. By definition of $v(a)$, this 
inequality is in fact an equality (almost everywhere), which concludes the 
proof.

Let us summarize now some properties of the volume:

\begin{prop} The function $a\to v(a)$ defined on the pseudoeffective cone is 
homogeneous of degree $n$, it satisfies $v(a+b)\geq v(a)+v(b)$ and it is 
continuous. Furthermore its restriction to the nef cone satisfies 
$v(a)=a^n$.
\end{prop}
Homogeneity and superadditivity being trivial, only the continuity needs a 
proof. But a consequence of the previous proof is the second of the 
inequalities $v(a)\leq\lim_{\ep\to 0}v(a+\ep\omega)\leq v(a)$, which easily 
implies the result.

\subsection{The Grauert-Riemenschneider criterion.}

We have seen in the first section that for a pseudoeffective line bundle $L$ 
on a compact K\"ahler manifold, the existence of a positive curvature 
current $T$ in $c_1(L)$ with $\int_XT_{ac}^n>0$ implies that $L$ is big. 
This is a Grauert-Riemenschneider-type criterion for bigness, which we would 
like to extend to any (i.e.~not necessarily rational) pseudoeffective class. 
As a matter of fact, let us say that a $(dd^c)$-cohomology $(1,1)$-class is 
$big$ if it contains a K\"ahler current. Then we prove the

\begin{theo}[Grauert-Riemenschneider criterion, K\"ahler case] If $a$ is a 
pseudoeffective $(1,1)$-class on the compact K\"ahler manifold $X$, then $a$ 
is big if, and only if, $v(a)>0$.
\end{theo}
This has been proven for a nef class by Demailly and Paun [DP01] as an 
essential step for their Nakai-Moishezon criterion for K\"ahler classes. 
Namely, they proved the

\begin{lem} If $a$ is a nef class $X$ with $\int_Xa^n>0$, then for every 
irreducible analytic subset $Y$ of $X$ of codimension $p$, there is some 
$\de>0$ such that $a^p\geq\de[Y]$, with $\geq$ meaning ``more 
pseudoeffective than''.
\end{lem}
Once this is proven, the following quite tricky argument also due to [DP01] 
enables us to conclude.

\noindent  Denote by $p_1$, $p_2$ the two projections $\ti X:=X\times X\to 
X$, and by $\Delta$ the diagonal. If $a$ is nef on $X$, then $\ti 
a=p_1^{\star}a+p_2^{\star}a$ is nef on $\ti X$, and it is easy to see that 
$\int_{\ti X}\ti a^{2n}={2n \choose n}(\int_Xa^n)^2>0$, so that by the above 
proposition one has for some $\de>0$:
$\ti a^n\geq\de[\Delta]$.

\noindent  But then for $\omega$ a K\"ahler form on $X$ we get:
$p_{1,\star}(\ti a^n\wedge p_2^{\star}\omega)\geq\de 
p_{1,\star}([\Delta]\wedge p_2^{\star}\omega)$.

\noindent
The first term in the inequality is easily seen to be 
$n(\int_Xa^{n-1}\wedge\omega)a$, whereas the second is $\de\omega$, whence 
we eventually get $a\geq(\de/(n\int_Xa^{n-1}\wedge\omega))\omega$, which is 
the expected conclusion.

What we plan to do is to apply the above strategy to the ``nef parts'' of 
``approximate Zariski decompositions'' of $a$. More concretely, reasoning as 
above, we can choose a sequence $T_k\in a$ of currents with analytic 
singularities such that $T_k\geq-\ep_k\omega$ and also that 
$v(T_{k,ac}+\ep_k\omega)$ is bounded away from zero. Then we consider finite 
sequences of blowing-ups $\mu_k:X_k\to X$ such that 
$\mu_k^{\star}T_k=\alpha_k+D_k$ with $\alpha_k$ smooth and $D_k$ a real 
effective divisor. Then the class 
$b_k:=\{\alpha_k+\ep_k\mu_k^{\star}\omega\}$ is nef and has volume bounded 
away from zero. Denote by $\ti X_k$ the product $X_k\times X_k$, by $p_1$ 
and $p_2$ the two projections, and $\Delta_k$ the diagonal. Consider 
moreover the class $\ti b_k:=p_1^{\star}b_k+p_2^{\star}b_k$. Then $\ti
b_k$ is nef, and has volume $v(\ti b_k)={2n \choose n}v(b_k)^2$ uniformly
bounded 
away from zero.\\ We claim the following

\begin{lem} There exists $\de>0$ such that $\ti b_k^n\geq\de[\Delta_k]$ for 
every $k$.
\end{lem}
Assuming for the moment that the lemma holds true, choose by means of Lemma 2.9 a sequence of smooth forms $u_k$ cohomologous to an effective linear combination $E_k$ of exceptional divisors of $\mu_k$uch that $\omega_k:=\mu_k^{\star}\omega-\eta_ku_k$ is a K\"ahler form 
for every $k$. Then Lemma 4.10 implies that $p_{1,\star}(\ti b_k^n\wedge 
p_2^{\star}\omega_k)\geq p_{1,\star}(\de[\Delta_k]\wedge 
p_2^{\star}\omega_k)$ on $X_k$. But the left hand side is 
$n(\int_{X_k}b_k^{n-1}\wedge\omega_k)b_k\leq 
n(\int_{X_k}b_k^{n-1}\wedge\omega_k)\mu_k^{\star}(a+\ep_k\omega)$ and the 
right hand side is $\de\omega_k$, thus pushing this forward by $\mu_k$, we 
eventually get 
$$n(\int_{X_k}b_k^{n-1}\wedge\omega_k)(a+\ep_k\omega)\geq\de\omega$$ on $X$. 
It remains to notice that 
$$\int_{X_k}b_k^{n-1}\wedge\omega_k=\{b_k\}^{n-1}\{\mu_k^{\star}\omega-[E_k]\}$$
$$\leq\int_Xb_k^{n-1}\wedge\mu_k^{\star}\omega=\int_X(\mu_{k,\star}\alpha_k+\ep_k\omega)^{n-1}\wedge\omega$$ 
(the inequality holds for $b_k$ is nef); now the class
$\{\mu_{k,\star}\alpha_k+\ep_k\omega\}$ is bounded, thus we have a uniform bound
for the integrals $\int_{X_k}b_k^{n-1}\wedge\omega_k$, and this is enough to
conclude.\\
It remains to prove the above lemma. To do this we use the following results
from [DP01]:

\begin{lem}[Concentration of the mass] Given a K\"ahler form $\omega$ on $X$ 
and $Y$ an analytic subset, there exists an almost plurisubharmonic function 
$\varphi$ on $X$ and a sequence $\varphi_{\ep}$ of smooth functions 
decreasing pointwise to $\varphi$ such that if we set 
$\omega_{\ep}:=\omega+dd^c\varphi_{\ep}$ then:

\noindent $(i)$ $\varphi$ has analytic singularities along $Y$,

\noindent $(ii)$ $\omega_{\ep}\geq\frac{1}{2}\omega$,

\noindent $(iii)$ For any smooth point $x\in Y$ at which $Y$ is 
$p$-codimensional and every neighbourhood $U$ of $x$, there exists $\de_U>0$ 
such that $\int_{U\cap V_{\ep}}\omega_{\ep}^p\wedge\omega^{n-p}\geq\de_U$ 
for every $\ep$, if $V_{\ep}=\{\varphi<\log\ep\}$.

\end{lem}

\begin{lem} Let $(X,\omega)$ be compact K\"ahler. Given an analytic subset 
$Y$, a function $\varphi$ and a family of closed smooth forms $\omega_{\ep}$ 
cohomologous to $\omega$ such that $(i)$, $(ii)$ and $(iii)$ in the above 
lemma hold, then for each nef class $a$ on $X$ such that $v(a)>0$, there 
exists $T\in a^p$ a closed positive current with
$$\int_{U\cap Y}T\wedge\omega^{n-p}\geq\eta$$
where $\eta=Cv(a)\de^2/Mv(\omega)>0$ with $C$ a universal constant,
$M=\int_Xa^{n-p}\wedge\omega^p$ and $\de:=\de_U$.
\end{lem}
The proof of 4.12 goes as follows: for each $\ep>0$, there exists by the Calabi-Yau 
theorem a smooth K\"ahler form $\alpha_{\ep}$ in the class $a+\ep\omega$ 
such that
$$\alpha_{\ep}^n=\frac{v(a+\ep\omega)}{v(\omega)}\omega_{\ep}^n.$$Denote by 
$\lambda^{\ep}_1\leq...\leq\lambda^{\ep}_n$ the eigenvalues of 
$\alpha_{\ep}$ with respect to $\omega_{\ep}$. Then we find:

\noindent$(a)$ $\lambda^{\ep}_1...\lambda^{\ep}_n=v(a+\ep\omega)/v(\omega)$.

\noindent$(b)$ 
$\alpha_{\ep}^p\geq\lambda^{\ep}_1...\lambda^{\ep}_p\omega_{\ep}^p$.

\noindent$(c)$ 
$\alpha_{\ep}^{n-p}\wedge\omega_{\ep}^p\geq\frac{p!(n-p)!}{n!}\lambda^{\ep}_{p+1
}...\lambda^{\ep}_n\omega_{\ep}^n$.\\By $(c)$, we get 
$$\int_X\lambda^{\ep}_{p+1}...\lambda^{\ep}_n\omega_{\ep}^n\leq {n\choose 
p}\{a+\ep\omega\}^{n-p}\{\omega\}^p=:M_{\ep},$$ thus in particular the set 
$E_{\eta}:=\{\lambda^{\ep}_{p+1}...\lambda^{\ep}_n\geq M_{\ep}/\eta\}$ has 
$\int_{E_{\eta}}\omega_{\ep}^n\leq\eta$ for $\eta>0$ small enough. Now
$(a)$ and $(b)$ yield 
$$\int_{U\cap V_{\ep}}\alpha_{\ep}^p\wedge\omega^{n-p}\geq\frac{v(a+\ep\omega)}{v(\omega)}\int_{U\cap V_{\ep}}\frac{1}{\lambda^{\ep}_{p+1}...\lambda^{\ep}_n}\omega_{\ep}^p\wedge\omega^{n-p}$$
$$\geq\frac{v(a+\ep\omega)}{v(\omega)}\int_{U\cap V_{\ep}-E_{\eta}}\frac{\eta}{M_{\ep}}\omega_{\ep}^p\wedge\omega^{n-p}.$$Observe that $$\int_{U\cap
  V_{\ep}-E_{\eta}}\omega_{\ep}^p\wedge\omega^{n-p}\geq\int_{U\cap V_{\ep}}\omega_{\ep}^p\wedge\omega^{n-p}-\int_{E_{\eta}}\omega_{\ep}^p\wedge\omega^{n-p}.$$The first integral is greater than $\de$ by assumption $(iii)$, and since
\\$\omega_{\ep}^p\wedge\omega^{n-p}\leq 2^{n-p}\omega_{\ep}^n$ by 
$(ii)$, the second integral will be less than $2^{n-p}\eta$. Combining all 
this yields in the end $$\int_{U\cap V_{\ep}}\alpha_{\ep}^p\wedge\omega^{n-p}\geq\frac{v(a+\ep\omega)}{v(\omega)}\frac{\eta}{M_{\ep}}(\de-2^{n-p}\eta).$$Now
we take $\eta:=\de/2^{n-p+1}$ and we choose $T$ some accumulation point of
$\alpha_{\ep}^p$. Then $T$ is a closed positive current in the cohomology class
$a^p$ such that $$\int_{U\cap Y}T\wedge\omega^{n-p}\geq C\frac{v(a)}{v(\omega)}\frac{\de^2}{M},$$ which concludes the proof.

We now denote by $\ti\mu_k:\ti X_k\to \ti X$ the product map 
$\mu_k\times\mu_k$, and we select on each $\ti X_k$ some K\"ahler form 
$\ti\omega_k$. First we apply Lemma 4.11 to $Y=\Delta$ in $\ti X$ and some 
neighbourhood $U$ of some point $x\in\Delta$ chosen to be a regular value of 
$\ti\mu_k$ for each $k$, so that we get $\varphi$, $\omega_{\ep}$ and 
$\de=\de_U$ as in the lemma. Then we consider for every $\rho>0$ the 
K\"ahler forms 
$\omega_{k,\ep}:=\ti\mu_k^{\star}\omega_{\ep}+\rho\ti\omega_k$ and 
$\omega_k:=\ti\mu_k^{\star}\omega+\rho\ti\omega_k$ and the function 
$\varphi_k:=\ti\mu_k^{\star}\varphi$ on $\ti X_k$. Then $(i)$, $(ii)$ and 
$(iii)$ still hold true, so that by Lemma 4.12 we get a closed positive $T_k\in 
\ti b_k^n$ such that

$$\int_{U_k\cap Y_k}T_k\wedge\omega_k^n\geq\eta_k$$ with 
$U_k:=\ti\mu_k^{-1}U$, $Y_k:=\ti\mu_k^{-1}Y$ and 
$\eta_k=Cv(b_k)\de^2/M_kv(\omega_k)$ for \\$M_k=\int_{\ti X_k}\ti 
b_k^n\wedge\omega_k^p$.
All this was depending on $\rho$, so we can now let $\rho$ tend to zero ($k$ 
is fixed) so as to get $\Theta_k$ some accumulation point of the 
$(T_k(\rho))_{\rho}$ with

$$\int_{U_k\cap Y_k}\Theta_k\wedge\ti\mu_k^{\star}\omega^n\geq\eta_k$$for 
$\eta_k=Cv(b_k)\de^2/M_kv(\ti\mu_k^{\star}\omega)$ and $M_k=\int_{\ti 
X_k}\ti b_k^n\wedge\ti\mu_k^{\star}\omega^n$. Now $M_k$ and 
$v(\ti\mu_k^{\star}\omega)=v(\omega)$ are certainly under uniform control, 
thus so is $\eta_k$, and the above equality compels $\Theta_k$ to have 
positive mass bounded away from zero on $\Delta_k$, since any other 
irreducible component of $Y_k$ is exceptional with respect to $\ti\mu_k$. 
Since $\Theta_k$ has bidimension $(n,n)$ and $\Delta_k$ is $n$-dimensional, 
Siu's decomposition formula then shows that $\Theta_k\geq\eta[\Delta_k]$ for 
some uniform $\eta$, which at last concludes the proof.

Let us now state some immediate corollaries of this bigness criterion:

\begin{prop} Let $f:X\to Y$ be a generically finite holomorphic map between 
compact K\"ahler manifolds, and $a$, $b$ be $(1,1)$-cohomology classes on 
$X$ and $Y$ respectively. Then one has:

\noindent$(i)$ $a$ big implies $f_{\star}a$ big,

\noindent$(ii)$ $b$ big implies $f^{\star}b$ big.

\end{prop}
The first point stems from the fact that pushing forward a K\"ahler current 
yields a K\"ahler current; as to the second, it is a consequence of 
Proposition 4.1 and the Grauert-Riemenschneider criterion.

\begin{prop} If $a$ is a pseudoeffective class with $v(a)>0$, there exists a 
sequence $T_k\in a$ of K\"ahler currents with analytic singularities such 
that $\int_XT_{k,ac}^n\to v(a)$ as $k\to+\infty$.
\end{prop}
Indeed, choose a positive $T\in a$ with $\int_XT_{ac}^n$ close to $v(a)$, and
$T_1\in 
a$ a K\"ahler current. Then for each $\ep>0$ small $T_{\ep}:=\ep 
T_1+(1-\ep)T$ is a K\"ahler current in $a$ with $T_{\ep,ac}=\ep 
T_{1,ac}+(1-\ep)T_{ac}$, thus $\int_XT_{\ep,ac}^n$ can be made as close to 
$v(a)$ as desired. Now we can approach each $T_{\ep}$ by K\"ahler currents 
with analytic singularities such that the absolutely continuous parts 
converge a.e., and we conclude by Fatou's lemma.

From the latter proposition, we deduce a general form of Fujita's theorem:

\begin{theo}[Approximation of the volume] Let $a$ be a pseudoeffective 
$(1,1)$-class on $X$ compact K\"ahler. Assume that the volume $v(a)$ is 
positive. Then for each $\ep>0$ small there exists a finite sequence of 
blowing-ups $\mu:\ti X\to X$, a K\"ahler class $\omega$ and a real effective 
divisor $D$ on $\ti X$ such that:

\noindent$(i)$ $\mu^{\star}a=\omega+D$

\noindent$(ii)$ $v(\omega)\leq v(a)\leq v(\omega)+\ep$

\noindent If $a$ is rational, $\omega$ and $D$ can be taken to be rational.
\end{theo}
$Remark$: of course the rational case is the original result of Fujita.

\noindent Proof of the theorem: select first $T\in a$ a K\"ahler current 
with analytic singularities such that $\int_XT_{ac}^n>v(a)-\ep$, as the 
proposition above allows; then take a resolution of singularities $\mu:\ti 
X\to X$ of $T$, so that $\mu^{\star}T=\alpha+D$ with $\alpha$ smooth and $D$ 
an effective divisor. If we denote by $E$ the exceptional divisor of $\mu$, 
Lemma 2.9 implies that the class of $\alpha-\de E$ is K\"ahler for each $\de>0$ 
small enough, thus the decomposition

$$\mu^{\star}a=(\alpha-\de E)+(D+\de E)$$
is the one sought after for $\de>0$ small enough, since $v(\alpha-\de 
E)\geq\int_{\ti X}\alpha^n=\int_XT_{ac}^n$. When $a$ is rational, on can 
arrange for $D$ to be rational by a slight perturbation which can be 
absorbed by $\alpha$, and $\de$ can be taken rational as well, q.e.d.

\subsection{Miscellaneous}

\subsubsection{The non-K\"ahler case}

Thanks to Lemma 4.2, the property of uniform boundedness of the mass is invariant
under modifications, thus still holds true on any compact complex manifold in
the Fujiki class $\cal C$ (a Fujiki manifold for short). As a consequence, all
the results which are not an equality of Calabi-Yau type remain true when the
K\"ahler assumption on $X$ is replaced by: $X$ is Fujiki. Since the presence of
a big class on $X$ forces it to be Fujiki by Proposition 2.8, we see for instance
that the general form of Fujita's theorem (Theorem 4.15) is true on an arbitrary
compact complex manifold if we assume $a$ to be big instead of having positive
volume. Nevertheless, it is tempting to think that the Grauert-Riemenschneider
criterion holds true in general, which comes down to the following conjecture
(cf.~also [DP01]):

\noindent{\bf Conjecture}: If a compact complex manifold $X$ carries a closed
positive $(1,1)$-current $T$ with $\int_XT_{ac}^n>0$, then $X$ is Fujiki. 

This conjecture is true for $\dim X=2$, as we will see.
\subsubsection{The case of a surface}

In this section, we assume that $X$ is a compact complex surface. Consider $T$ a
closed positive $(1,1)$-current on $X$, and let $T=D+R$ stand for its Siu
decomposition. Then we claim that the $dd^c$-class $\{R\}$ is nef, and is even
K\"ahler when $T$ is a K\"ahler current. This is in fact a straightforward
consequence of a result of Paun given in [DP01], which says that a
pseudoeffective class (resp.~a big class) $\{T\}$ is nef (resp.~K\"ahler) if its
restriction $\{T\}_{|Y}$ is nef (resp.~K\"ahler) for every irreducible analytic
subset contained in some Lelong sublevel set $E_c(T)$. We will however briefly
recall the argument: let $R_k$ be an approximation of $R$ with analytic
singularities as in Theorem 2.6. Then $R_k$ will be smooth outside a finite
subset, and we have $R_k\geq-\ep_k\omega$ (resp.~$R_k\geq\ep\omega$ when $R$ is
a K\"ahler current). Near each singular point of $R_k$, let us write
$R_k=dd^c\varphi_k$ for some local potential $\varphi_k$. Then we replace
$\varphi_k$ by $\max_{\eta}(\varphi_k,-C_k)$, with $C_k>0$ big enough so that
the gluing condition is satisfied, and the new $R_k$ thus obtained is smooth
with $R_k\geq-\ep_k\omega$. In the second case where $R_k\geq\ep\omega$, we
choose some local potential $\varphi_k$ near each singular point, and replace it
by $\max_{\eta}(\varphi_k,|z|^2-C_k)$ with $C_k>0$ big enough so that the gluing
condition is satisfied, and we get in this way a smooth $R_k$ with a strictly
positive lower bound, i.e.~a K\"ahler form. 
The above reasoning shows in particular that $X$ is Fujiki iff it is K\"ahler.

Now if we assume that $\int_X T_{ac}^2>0$, then since $\{R\}$ is nef we get by
Lemma~4.4 $\{R\}^2\geq\int_XR_{ac}^2=\int_XT_{ac}^2>0$, thus the intersection form
on $H^{1,1}(X,\er)$ cannot be negative definite. Now this compels $b_1(X)$ to be
even by classical results of Kodaira (cf.~e.g.~[Lam99]), and this in turn
forces $X$ to be K\"ahler by the main result of [Lam99] or [Buc99]. Therefore
the above conjecture is true on a surface. In particular, whenever $v(a)>0$ for
some pseudoeffective class $a$, we can choose $T\in a$ such that
$\int_XT_{ac}^2=v(a)$, and then we get $v(a)\geq
v(\{R\})=\{R\}^2\geq\int_XT_{ac}^2=v(a)$, thus $a$ and $\{R\}$ have same volume.
Let us summarize all this in the

\begin{prop} Let $X$ be a compact complex surface. Then for each pseudoeffective
class $a$ on $X$, there exists a ``Zariski decomposition'' $a=d+r$ where $d$ is
an infinite series of real effective divisors and $r$ is a nef class with $v(a)=r^2$.
Furthermore, as soon as $v(a)>0$, $a$ is big and $r$ can be taken to be a
K\"ahler class.
\end{prop}

\subsubsection{Behaviour of the volume in deformations}

In this last part we prove the following

\begin{prop} The volume is upper-semicontinuous on K\"ahler
  deformations in the following sense: if $\mathcal{X}\to S$ is a proper
  submersive K\"ahler map and $a$ is a pseudoeffective $(1,1)$-class
  on the central fibre $X_0=:X$, then one has $$v(a)\geq\limsup_{b\to
a}v(b)$$where $b$ is a pseudoeffective $(1,1)$-class on some fibre $X_t$.
\end{prop}
$Proof$: Upon shrinking the base $S$, we may assume that the deformation is
topologically trivial, and that there exist K\"ahler metrics $\omega_t$ on $X_t$
depending smoothly on $t$; we normalize them so that $\int_{X_t}\omega_t^n=1$.
Let then $b_k$ be a sequence of cohomology classes in $H^2(X,\co)$ converging to
$a$, and such that $b_k$ is a pseudoeffective $(1,1)$-class on $X_{t_k}$ for
some sequence $t_k$ going to 0, and that $v(b_k)$ converges to the right hand
side of our inequality. By the degenerate Calabi-Yau theorem, we can choose a
closed positive $(1,1)$-current $T_k\in b_k$ on $X_{t_k}$ such that
$T_{k,ac}^n=v(b_k)\omega_{t_k}^n$ almost everywhere. Since $b_k$ converges, it
is bounded, thus $T_k$ is bounded in mass, and after extracting some
subsequence, we may assume that $T_k$ converges to a closed positive $T\in a$.
Then we apply Lemma 3.5 to get in the limit $T_{ac}^n\geq\limsup_{b\to
a}v(b)\omega^n$, which
gives the expected result. 

As a consequence, we get the following ``quantitative'' version of a
result of Huybrechts [Huy01]: denote by $\mathcal{C}_X$ the connected
component of the K\"ahler cone in the open cone $\{a\in
H^{1,1}(X,\er),\int_Xa^n>0\}$. Then we have:

\begin{prop} If $X$ is an irreducible symplectic (K\"ahler) $2n$-dimensional
manifold, then $v(a)\geq\int_Xa^{2n}$ for each $a\in\mathcal{C}_X$.
\end{prop}  

Indeed, we can choose a sequence of very general points $t_k$ converging to 0 in
the universal deformation space of $X$ such that $a$ is a limit of
$(1,1)$-classes $a_k$ on $X_{t_k}$. By [Huy01], the cone $\mathcal{C}_{X_{t_k}}$
coincides with the K\"ahler cone of $X_{t_k}$ because $t_k$ is very general,
therefore $a_k$ is a K\"ahler class and we have $v(a_k)=\int_Xa_k^{2n}$ for each
$k$. The result then follows from Proposition 4.17.   

\it{I would like here to address my warmest thanks to my thesis advisor
Jean-Pierre Demailly for his many suggestions and his interest in this work.}

\section{References.}
\begin{itemize}
\item{\bf [Bon93]} Bonavero, L.\ --- {\sl In\'egalit\'es de Morse
holomorphes singuli\`eres}, C.\ R.\ Acad.\ Sci.\ S\'erie I {\bf 317} (1993),
1163--1166.

\item{\bf [Buc99]} Buchdahl, N.\ --- {\sl On compact K\"ahler surfaces}, Ann.
Inst. Fourier {\bf 50} (1999), 287--302. 

\item{\bf [Dem82]} Demailly, J.-P.\ --- {\sl Estimations $L^2$ pour
l'op\'erateur $\overline{\partial}$ d'un fibr\'e vectoriel holomorphe
semi-positif au dessus d'une vari\'et\'e k\"ahlerienne compl\`ete}, Ann. Sci.
Ecole Norm. Sup. {\bf 15} (1982), 457--511.

\item{\bf [Dem92]} Demailly, J.-P.\ --- {\sl Regularization of closed positive
currents and intersection theory}, J. Alg. Geom. {\bf 1} (1992), 361--409.

\item{\bf [Dem93]} Demailly, J.-P.\ --- {\sl A numerical criterion for very
ample line bundles}, J. Diff. Geom. {\bf 37} (1992), 323--374.

\item{\bf [Dem96]} Demailly, J.-P. \ --- {\sl $in$ Introduction \`a la th\'eorie
de Hodge}, Panoramas et synth\`eses, S.M.F. {\bf 3} (1996), 3--111.

\item{\bf [DEL00]} Demailly, J.-P.; Ein, L.; Lazarsfeld, R.\ --- {\sl A
subadditivity property of multiplier ideals}, math.AG/0002035 (2000). 

\item{\bf [DP01]} Demailly, J.-P.; Paun, M.\ --- {\sl Numerical characterization
of the K\"ahler cone of a compact K\"ahler manifold}, math.AG/0105176 (2001).

\item{\bf [Fuj94]} Fujita, T.\ --- {\sl Approximating Zariski decomposition of
big line bundles}, Kodai Math. J. {\bf 17} (1994), 1--3.

\item{\bf [Huy01]} Huybrechts, D.\ --- {\sl $Erratum$: Compact hyperk\"ahler
manifolds: basic results}, math.AG/0106014 (2001).

\item{\bf [JS93]} Ji, S.; Shiffman, B.\ --- {\sl Properties of compact complex
manifolds carrying closed positive currents}, J. Geom. Anal. {\bf 3}, 1 (1993),
37--61.

\item{\bf [Lam99]} Lamari, A.\ --- {\sl Courants k\"ahleriens et surfaces
compactes}, Ann. Inst. Fourier {\bf 49} (1999), 249--263.

\item{\bf [Mou98]} Mourougane, C.\ --- {\sl Versions k\"ahleriennes du
th\'eor\`eme d'annulation de Bogomolov}, Collect. Math. {\bf 49}, 2-3 (1998),
433--445.

\item{\bf [Yau78]} Yau, S.-T.\ --- {\sl On the Ricci curvature of a 
complex K\"ahler manifold and the complex Monge--Amp\`ere equation}, Comm. Pure
Appl. Math. {\bf 31} (1978).
\end{itemize}
\end{document}